# STABILITY IN DISTRIBUTION OF RANDOMLY PERTURBED QUADRATIC MAPS AS MARKOV PROCESSES

By Rabi Bhattacharya[1] and Mukul Majumdar

*Indiana University and Cornell University*

Iteration of randomly chosen quadratic maps defines a Markov process: $X_{n+1} = \varepsilon_{n+1} X_n (1 - X_n)$, where $\varepsilon_n$ are i.i.d. with values in the parameter space $[0, 4]$ of quadratic maps $F_\theta(x) = \theta x(1 - x)$. Its study is of significance as an important Markov model, with applications to problems of optimization under uncertainty arising in economics. In this article a broad criterion is established for positive Harris recurrence of $X_n$.

**1. Introduction.** The present article explores the problem of stability in distribution of randomly perturbed dynamical systems of random compositions of the form

$$(1.1) \qquad X_n = F_{\varepsilon_n} F_{\varepsilon_{n-1}} \cdots F_{\varepsilon_1} X_0, \qquad n \geq 1,$$

where $\varepsilon_n$, $n \geq 1$, is an i.i.d. sequence with values in the parameter space $[0, 4]$ of the quadratic map $F_\theta(x) = \theta x(1 - x)$, and $X_0$ is independent of $\{\varepsilon_n : n \geq 1\}$. To exclude the trivial invariant probability $\delta_{\{0\}}$, we will restrict the state space of the Markov process to $S = (0, 1)$. By *stability in distribution* we mean the convergence in distribution of $\frac{1}{n} \sum_{m=0}^{n-1} p^m(x, dy)$ to the same limit $\pi(dy)$ for every initial state $x \in S$, with $p^{(n)}$ denoting the $n$-step transition probability of $\{X_n : n \geq 0\}$. Then $\pi$ is the unique invariant probability of this Markov process.

One may, more generally, consider a stationary ergodic sequence $\{\varepsilon_n : n \geq 1\}$, rather than an i.i.d. one. For some results in this direction we refer to Anantharam and Konstantopoulos (1997). Our motivations for using the present framework are briefly stated below.

We mention one particular application from economics that has provided at least a part of the motivation for this work. Consider a dynamic optimization problem in which one is given a "*production function*" $f : \mathbb{R}_+ \to \mathbb{R}_+$

Received December 2002; revised September 2003.
[1]Supported in part by NSF Grant DMS-00-73865 and a Guggenheim Fellowship.

and a *welfare function* $w:\mathbb{R}_+^2 \times A \to \mathbb{R}_+$, where $A$ is a parameter set, say $A = [1,4]$, parametrizing a family of economies. For an initial $x \geq 0$, a *program* $x_n : n \geq 0$, is a sequence such that $0 \leq x_0 = x \leq x_n \leq f(x_{n-1})$. The *consumption sequence* $\{c_n : n \geq 1\}$ is defined as $c_n = f(x_{n-1}) - x_n$. Given a *discount factor* $\delta > 0$ and a parameter value $\theta$, one wishes to find an *optimal program* $\{\hat{x}_n : n \geq 0\}$, $\hat{x}_0 = x$, which maximizes $\sum_{n=0}^\infty \delta^n w(x_n, c_{n+1}, \theta)$ over all programs $\{x_n : n \geq 0\}$ starting at $x_0 = x$. Under economically feasible assumptions one may find $f$ and $w$ such that an optimal program is given recursively by $\hat{x}_{n+1} = \theta \hat{x}_n (1 - \hat{x}_n), x \in [0,1]$ [Majumdar and Mitra (2000); also see Bala and Majumdar (1992) for a special case]. That is, the optimal program is given by the trajectory of the dynamical system $F_\theta$ with initial state $x$. Since "uncertainty" is inherent in economic systems, one may thus obtain a randomly perturbed quadratic system. Alternatively, one may at the outset consider a *stochastic dynamic programming problem* and directly arrive at a stationary optimal policy leading to an evolution of states of the form (1.1) [Mitra (1998)].

For earlier work on the existence of a unique invariant probability of the process (1.1) when the distribution of $\varepsilon_n$ has a two-point support or, more generally, may not have a density component, see Bhattacharya and Majumdar (1999), Bhattacharya and Rao (1993), Bhattacharya and Waymire (2002) and Carlsson (2002). Since $\{X_n\}$ is then generally non-Harris or nonirreducible, fairly strong restrictions on the distribution of $\varepsilon_n$ are needed to ensure uniqueness of the invariant probability. The question of Harris recurrence of the process (1.1) was raised by Athreya and Dai (2000).

For ease of reference, we recall that a Markov process on $(S, \mathbb{S})$, or its transition probability $p(x, dy)$, is *irreducible* if there exists a nonzero $\sigma$-finite measure $\varphi$ such that

$$(1.2) \qquad U(x, B) := \sum_{n=1}^\infty p^{(n)}(x, B) > 0 \qquad \forall\, x \in S \text{ if } \varphi(B) > 0.$$

In this case one also says that the Markov process is $\varphi$-*irreducible*. A $\varphi$-irreducible Markov process is *Harris recurrent*, or $\varphi$-*recurrent*, if

$$(1.3) \qquad U(x, B) = \infty \qquad \forall\, x \in S \text{ if } \varphi(B) > 0.$$

A $\varphi$-irreducible Markov process is said to be *Harris positive* if it has an invariant probability $\pi$. It is known that a Harris positive process is Harris recurrent [Meyn and Tweedie (1993), Proposition 10.1.1]. Therefore, one may call the process *Harris positive recurrent* in this case, or *Harris ergodic*. The invariant probability measure $\pi$ of a Harris positive (recurrent) process is unique and one has the convergence of the Caesaro mean in total variation distance,

$$(1.4) \qquad \sup \left| \frac{1}{n} \sum_{m=1}^n p^{(m)}(x, B) - \pi(B) \right| \to 0 \qquad \text{as } n \to \infty \ \forall\, x \in S.$$

In our present content irreducibility is proved with the help of Lemmas 2.3 and 2.4. A tightness assumption does the rest. A verifiable sufficient condition for tightness due to Athreya and Dai (2000) then yields the useful Corollary 2.2.

Before concluding this section, let us mention that our goal is to derive Harris positive recurrence of the process (1.1) under as broad a (verifiable) condition as we can muster. That this task is rather delicate is perhaps evident from Remark 2.4. In addition to the economic application mentioned above, another motivation for this work comes from a program proposed originally by Kolmogorov for the approximate computation of the physically meaningful ergodic invariant probability of a dynamical system. Note that a chaotic dynamical system has infinitely many extremal, or ergodic, invariant probabilities (including one uniform distribution on each of its infinitely many periodic orbits). The problem is to approximate the invariant probability to which most points are attracted. Kolmogorov had proposed that one introduce a small random (possibly absolutely continuous) perturbation to the dynamical system, such that the resulting Markov process has a unique invariant probability, which is an approximation to the so-called *Kolmogorov measure* of the dynamical system [see Kifer (1988) and Katok and Kifer (1986)]. In our context this may be achieved, for example, by having a uniform density on a small interval around a parameter point of interest. Extensive simulations have shown that, for the present case, this program works not only for the chaotic regime, but also for $\theta$ such that $F_\theta$ has a stable periodic orbit, or a quasi-periodic attractor. These simulations may be found in Bhattacharya and Majumdar (2002).

**2. Harris recurrence and ergodicity.** On the state space $S = (0, 1)$, consider the Markov process defined recursively by

$$(2.1) \qquad X_{n+1} = F_{\varepsilon_{n+1}} X_n, \qquad n = 0, 1, 2, \ldots,$$

where $\{\varepsilon_n : n \geq 1\}$ is a sequence of i.i.d. random variables with values in $(0, 4)$ and, for each value $\theta \in (0, 4)$, $F_\theta$ is the quadratic function (on $S$):

$$(2.2) \qquad F_\theta x \equiv F_\theta(x) = \theta x(1 - x), \qquad 0 < x < 1.$$

As always, the initial random variable $X_0$ is independent of $\{\varepsilon_n : n \geq 1\}$. Our main result provides a criterion for Harris recurrence and the existence of a unique invariant probability for the process $\{X_n : n \geq 0\}$. Recall that a sequence $\mu_n$, $n \geq 1$, of probability measures on $S$ is said to be *tight* if, for every $\varepsilon > 0$, there exists a compact $K_\varepsilon \subset S$ such that $\mu_n(K_\varepsilon) \geq 1 - \varepsilon$ for all $n \geq 1$.

Let $p(x, dy)$ denote the (one-step) transition probability of $\{X_n : n \geq 1\}$ and $p^{(m)}(x, dy)$ the corresponding $m$-step transition probability.

THEOREM 2.1. *Assume that the distribution of $\varepsilon_1$ has a nonzero absolutely continuous component [w.r.t. Lebesgue measure on $(0,4)$] whose density is bounded away from zero on some nondegenerate interval in $(1,4)$. If, in addition, $\{\frac{1}{N}\sum_{n=1}^{N} p^{(n)}(x,dy) : N \geq 1\}$ is tight on $S = (0,1)$ for some $x$, then:*

(i) *$\{X_n : n \geq 0\}$ is Harris recurrent and has a unique invariant probability $\pi$,*

(ii) *$\frac{1}{N}\sum_{n=1}^{N} p^{(n)}(x,dy)$ converges to $\pi$ in total variation distance, for every $x$, as $n \to \infty$.*

COROLLARY 2.2. *If $\varepsilon_1$ has a nonzero density component which is bounded away from zero on some nondegenerate interval contained in $(1,4)$ and if, in addition,*

$$(2.3) \qquad E\log\varepsilon_1 > 0 \quad \text{and} \quad E|\log(4-\varepsilon_1)| < \infty,$$

*then $\{X_n : n \geq 0\}$ has a unique invariant probability $\pi$ on $S = (0,1)$ and $(1/N)\sum_{n=1}^{N} p^{(n)}(x,dy) \to \pi$ in total variation distance, for every $x \in (0,1)$.*

Note that if the support of the distribution $Q$ of $\varepsilon_1$ is contained in $[\mu, \nu]$, where $1 < \mu < \nu < 4$, then $[a,b] \equiv [\min\{1 - \frac{1}{\mu}, F_\mu(\frac{\nu}{4})\}, \frac{\nu}{4}]$ is an invariant interval for the Markov process (2.1) [see Bhattacharya and Rao (1993) and Bhattacharya and Waymire (2002)]. Since the transition probability has the Feller property, whatever be $Q$, there exists an invariant probability with support contained in $[a,b] \subset (0,1)$. The result of Athreya and Dai (2000) is an important generalization of this.

We will need some lemmas for the proof of this theorem.

REMARK 2.1. Note that the tightness condition in Theorem 2.1, guaranteed, for example, by (2.3), cannot be dispensed with. It has been shown by Athreya and Dai (2000) that if $E\log\varepsilon_1 \leq 0$, then $X_n$ converges in probability to 0.

LEMMA 2.3. *Suppose the distribution $Q$ of $\varepsilon_1$ on $(0,4)$ has a nonzero absolutely continuous component (w.r.t. Lebesgue measure $\lambda$) whose density $h(\theta)$ is bounded away from zero on an interval $[c,d]$, $1 < c < d < 4$. Then there exist a nonempty open interval $J \subset (0,1)$, a number $\delta > 0$ and a positive integer $m$ such that*

$$(2.4) \qquad \inf_{x \in J} p^{(m)}(x, B) \geq \delta\lambda(B) \qquad \forall\ \text{Borel}\ B \subset J.$$

PROOF. First assume $Q$ is absolutely continuous with a continuous density $h$. Let $\theta_0 \in (1,4)$ be such that $h(\theta_0) > 0$ and $F_{\theta_0}$ has an attractive periodic orbit of period $m$. Such a point $\theta_0$ exists, since the set of points $\theta$ for which $F_\theta$ has an attractive fixed point or an attractive periodic orbit is dense in $(0, 4)$, by a result of Graczyk and Swiatek (1997). The $n$-step transition probability density is continuous in $(x, y)$ and is given recursively by

$$p(x,y) \equiv p^{(1)}(x,y) = \frac{1}{x(1-x)} h\left(\frac{y}{x(1-x)}\right),$$

(2.5) $$p^{(n+1)}(x,y) = \int_{(0,1)} \frac{1}{z(1-z)} h\left(\frac{y}{z(1-z)}\right) p^{(n)}(x,z)\, dz,$$

$$x, y \in (0,1), n \geq 1.$$

Let $\{x_0, x_1, \ldots, x_{m-1}\}$ be the attractive periodic orbit of $F_{\theta_0}: F_{\theta_0} x_{i-1} = x_i, i = 1, \ldots, m, x_m \equiv x_0$. Then

(2.6) $$p^{(1)}(x_{i-1}, x_i) = \frac{1}{x_{i-1}(1 - x_{i-1})} h\left(\frac{x_i}{x_{i-1}(1 - x_{i-1})}\right)$$
$$= \frac{h(\theta_0)}{x_{i-1}(1 - x_{i-1})} > 0, \qquad 1 \leq i \leq m,$$

since $x_i = \theta_0 x_{i-1}(1 - x_{i-1}) \equiv F_{\theta_0} x_{i-1}$. By (2.6) and the continuity of $(x, y) \to p^{(1)}(x, y)$, there exist $\delta_i > 0$ such that

$$g(y_1, \ldots, y_{m-1}) := p^{(1)}(x_0, y_1) p^{(1)}(y_1, y_2) \cdots p^{(1)}(y_{m-2}, y_{m-1}) p^{(1)}(y_{m-1}, x_0)$$
$$> 0 \qquad \forall\, y_i \in [x_i - \delta_i, x_i + \delta_i], 1 \leq i \leq m-1,$$

so that

(2.7) $$p^{(m)}(x_0, x_0) \geq \int \cdots \int g(y_1, \ldots, y_{m-1})\, dy_1 \cdots dy_{m-1} > 0,$$

where the integration is over the rectangle $[x_1 - \delta_1, x_1 + \delta_1] \times \cdots \times [x_{m-1} - \delta_{m-1}, x_{m-1} + \delta_{m-1}]$. By the continuity of $(x, y) \to p^{(m)}(x, y)$, it follows that there exists an open neighborhood $J$ of $x_0$ such that

(2.8) $$p^{(m)}(x, y) \geq \delta > 0 \qquad \forall\, x, y \in \bar{J},$$

where $\bar{J}$ is the closure of $J$ in $(0, 1)$. This proves (2.4) assuming that $Q$ is absolutely continuous with a continuous density. In the general, case let $I \subset (1, 4)$ be a nondegenerate closed interval such that $h(\theta) \geq \delta' > 0 \,\forall\, \theta \in I$. There exists a nonnegative continuous function $\underline{h}$ on $(0, 4)$ such that $\underline{h}(\theta) > 0 \,\forall\, \theta \in$ interior of $I$, and $\underline{h}(\theta) \leq h(\theta) \,\forall\, \theta \in (0, 4)$. Define $\underline{p}^{(n)}(x, y)$ in place of $p^{(n)}(x, y)$, $n \geq 1$, in (2.5) by replacing $h$ by $\underline{h}$. Let $\theta_0$ be a point in the interior of $I$ such that $F_{\theta_0}$ has an attractive periodic orbit of period $m$, say, $\{x_0, x_1, \ldots, x_{m-1}\}$. Then the same argument as given above shows that there exists an open neighborhood $J$ of $x_0$ such that $\underline{p}^{(m)}(x, y) \geq \delta \,\forall\, x, y \in \bar{J}$, for some $\delta > 0$. Since $h \geq \underline{h}$, $p^{(m)}(x,y) \geq \underline{p}^{(m)}(x,y) \geq \delta \ \forall x, y \in \bar{J}$, and the proof of (2.4) is complete. $\square$

Our final lemma adds greater specificity to Lemma 2.3 and to the proof of Theorem 2.1.

LEMMA 2.4. *Assume the hypothesis of Lemma 2.3. There exist $\gamma_1, \gamma_2$, $c < \gamma_1 < \gamma_2 < d$, and $m \geq 1$ such that:*

(a) *$F_\theta$ has an attractive periodic orbit of period $m$ for every $\theta \in (\gamma_1, \gamma_2)$, and*

(b) *if $q(\theta)$ denotes the largest point of the attractive periodic orbit of $F_\theta(\theta \in (\gamma_1, \gamma_2))$, then there exists an open interval $J \subset (0,1)$ for which:*

(i) (2.4) *holds,*
(ii) $q(\theta) \in J \ \forall \theta \in (\gamma_1, \gamma_2)$,
(iii) $\theta \to q(\theta)$ *is a diffeomorphism on $(\gamma_1, \gamma_2)$ onto $J$.*

PROOF. As in the proof of Lemma 2.3, let $\theta_0 \in (c,d)$ be such that $F_{\theta_0}$ has an attractive periodic orbit of some period, say, $m$. Apply the inverse function theorem to the function $(\theta, x) \to F_\theta^m x - x$ in a neighborhood of $(\theta_0, q(\theta_0))$. For this, note that

$$(2.9) \qquad \left(\frac{d}{dx}\{F_\theta^m x - x\}\right)_{\theta=\theta_0, x=q(\theta_0)} < 0,$$

in view of the property $|\frac{d}{dx}F_{\theta_0}^m x|_{x=q(\theta_0)} < 1$ [since $q(\theta_0)$ is an attractive fixed point of $F_{\theta_0}^m$]. Hence there exists $\underline{\theta} < \theta_0 < \bar{\theta}$ such that $\theta \to q(\theta)$ is a diffeomorphism on $(\underline{\theta}, \bar{\theta})$ onto an open interval $I \subset (0,1)$. Now apply Lemma 2.3 to find an open interval $J = (u_1, u_2) \subset I, u_1 < q(\theta_0) < u_2$, such that (2.4) holds, and let $\gamma_i = q^{-1}(u_i)$, $i = 1, 2$. $\square$

PROOF OF THEOREM 2.1. Let $\pi$ be an ergodic (i.e., extremal) invariant probability on $S = (0,1)$, which exists by the assumption of tightness. We will first show that $\pi(J) > 0$ for the set $J$ in Lemma 2.4. Fix $x \in (0,1)$. There exists a point in the interval $(F_{\gamma_1} x, F_{\gamma_2} x)$ which is attracted to the (attractive) periodic orbit of $F_{\theta_0}$, where $\theta_0$ is as in the proof of Lemma 2.4. Note that, outside a set of Lebesgue measure zero, every point of $(0,1)$ is so attracted [see, e.g., Collet and Eckmann (1980), page 13]. Thus there exist $n$ and $\theta_1^0, \theta_2^0, \ldots, \theta_n^0 \in (\gamma_1, \gamma_2)$ such that $F_{\theta_1^0} F_{\theta_2^0} \cdots F_{\theta_n^0} x \in J$. Consider the open subset of $(0,1) \times (\gamma_1, \gamma_2)^n$ given by $\{(y, \theta_1, \theta_2, \ldots, \theta_n) : F_{\theta_1} F_{\theta_2} \cdots F_{\theta_n} y \in J\}$. Since $(x, \theta_1^0, \theta_2^0, \ldots, \theta_n^0)$ belongs to this open set, there exists a neighborhood of this point, say, $(y_1, y_2) \times (\theta_{11}, \theta_{12}) \times \cdots \times (\theta_{n1}, \theta_{n2}) \subset (0,1) \times (\gamma_1, \gamma_2)^n$ such that $\forall (y, \theta_1, \ldots, \theta_n)$ in this neighborhood, $F_{\theta_1} F_{\theta_2} \cdots F_{\theta_n} y \in J$. This implies that, for every initial state $y \in (y_1, y_2) = I_x$, say, the probability that in $n$ steps the Markov process will reach $J$ is at least $c_1^n(\theta_{12} - \theta_{11}) \cdots (\theta_{n2} - \theta_{n1}) \equiv \varepsilon_n(x) > 0$, where $c_1 := \inf\{h(\theta) : \theta \in [c, d]\}$. Now choose $x$ such that it belongs to the support of $\pi$. Then $\pi(I_x) > 0$ and, with $n = n(x)$ as above,

$$\begin{aligned}(2.10) \quad \pi(J) &= \int p^{(n)}(y, J)\pi(dy) \geq \int_{I_x} p^{(n)}(y, J)\pi(dy) \\ &\geq \varepsilon_n(x)\pi(I_x) > 0.\end{aligned}$$

Also, by Lemma 2.3,

$$(2.11) \quad \pi(B) \geq \int_J p^{(m)}(x, B)\pi(dx) \geq \delta\lambda(B)\pi(J) \qquad \forall \text{Borel } B \subset J.$$

In particular, $\pi$ has an absolutely continuous component on $J$ w.r.t. Lebesgue measure, with a density bounded below by $\delta\pi(J) > 0$. Since the same argument would apply to every invariant ergodic probability $\pi_1$, and two distinct extremal invariant measures are mutually singular, it follows that $\pi$ is the unique invariant probability.

Since we have argued above that for every $x \in S = (0, 1)$ there exists $n = n(x)$ such that $p^{(n)}(x, J) > 0$, we get, using (2.4),

$$\begin{aligned}(2.12) \quad p^{(n+m)}(x, B) &\geq \int_J p^{(m)}(z, B)p^{(n)}(x, dz) \\ &\geq \delta\lambda(B)p^{(n)}(x, J) \\ &> 0 \qquad \forall x \in S \equiv (0, 1), \forall B \subset J, \lambda(B) > 0.\end{aligned}$$

Hence the Markov process is *irreducible* with respect to the measure $\phi(B) := \lambda(B \cap J), B$ Borel $\subset (0, 1)$. From the standard theory for Harris processes [see, e.g., Meyn and Tweedie (1993), Proposition 10.1.1, page 231], it now follows that the Markov process (2.1) is positive Harris recurrent. Part (ii) of Theorem 2.1 is a consequence of this fact [see, e.g., Meyn and Tweedie (1993), Theorem 10.1.1, page 230, and Orey (1971)]. □

REMARK 2.2. Under the hypothesis of Theorem 2.1, the Markov process is not in general aperiodic. For example, one may take the distribution of $\varepsilon_n$ to be concentrated in an interval such that, for every $\theta$ in this interval, $F_\theta$ has a stable periodic orbit of period $m > 1$. As may be seen from Lemma 2.4, one may find an interval of this kind so that the process is irreducible and cyclical of period $m$. If $\varepsilon_n$ has a density component bounded away from zero on a nondegenerate interval $B$ containing a *stable fixed point*, that is, $B \cap (0, 3) \neq \varnothing$, then the process is *aperiodic* and $p^{(n)}(x, \cdot)$ converges in total variation distance to a unique invariant $\pi$. Restrictive assumptions of this kind have been used by Bhattacharya and Rao (1993) and Dai (2000).

REMARK 2.3. In the course of the proof of Theorem 2.1, we showed that there exists an open interval $J$ which is recurrent. An alternative way to complete the proof, in presence of Lemma 2.3, would then be to prove that $J$ is *positive recurrent* in the sense that $\sup\{E\tau_J^x : x \in J\} < \infty$. Here $\tau_J^x := \inf\{n \geq 1 : X_n \in J | x_0 = x\}$. Under our general assumption, this last task appears elusive, compared to the route via irreducibility.

REMARK 2.4. We do not know if the conclusion of Theorem 2.1 remains valid under the assumption "the distribution $Q$ (of $\varepsilon_1$) on $(0, 4)$ has a nonzero absolutely continuous component with respect to Lebesgue measure on $(1, 4)$," in addition to (2.3). Note that such a $Q$ may assign its entire mass on the set of all $\theta$ for which $F_\theta$ is chaotic.

**Acknowledgment.** The authors wish to thank the referee for valuable suggestions.

| | |
|---|---|
| Department of Mathematics | Department of Economics |
| University of Arizona | Cornell University |
| 617 N. Santa Rita Avenue | 460 Uris Hall |
| Tucson, Arizona 85721 | Ithaka, New York 14853 |
| USA | USA |
| e-mail: rabi@math.arizona.edu | e-mail: mkm5@coruell.edu |